\documentclass[10pt]{amsart}
\usepackage{amssymb,latexsym,amsmath,amsfonts}

\numberwithin{equation}{section}
\theoremstyle{definition}
\newtheorem{definition}{Definition}[section]
\newtheorem{remark}[definition]{Remark}
\theoremstyle{plain}
\newtheorem{theorem}[definition]{Theorem}
\newtheorem{lemma}[definition]{Lemma}

\newcommand{\dia}{{\rm r}_x}
\newcommand{\Cn}{{\mathbb{C}}^n}
\newcommand{\CC}{{\mathbb{C}}^2}
\newcommand{\RN}{{\mathbb{R}}^n}

\newcommand{\OM}{\overline\Omega}
\newcommand{\unal}{\mathcal{E}}

\newcommand{\eps}{\varepsilon}
\newcommand{\sig}{\sigma}

\newcommand{\vaph}{1-\alpha}

\newcommand{\smoo}{\mathcal{C}}

\begin{document}
\title[Peak functions]{A procedure for constructing peak functions}
\author{Gautam Bharali}
\address{Department of Mathematics, The University of Michigan, 525 East University Avenue,
Ann Arbor, MI 48109}
\email{bharali@umich.edu}
\thanks{Supported in part by a Rackham Faculty Fellowship, The Rackham Graduate School,
University of Michigan}

\begin{abstract} We extend Bishop's one-fourth three-fourths principle for constructing
peak functions belonging to a uniform algebra to a situation where the ``approximate
barriers'' associated with the Bishop construction are not uniformly bounded.
\end{abstract}

\maketitle

\section{Introduction}\label{S:intro}

In this paper, we revisit a procedure for constructing peak functions devised by Bishop.
Specifically, let $\Omega$ be a bounded domain in
$\RN$ (or $\Cn$), and let $\unal$ be a closed subspace of $\smoo(\OM)$ : the class of
all complex-valued functions that are continuous on $\OM$.
Let $x\in\OM$; we say that $\boldsymbol{f}$ (here $f\in \unal$) {\bf peaks at} 
$\boldsymbol{x}$ if $f(x)=1$ and $|f(y)|<1 \ \forall y\in \OM\setminus\{x\}$,
and we call $f$ a {\bf peak function of class} $\boldsymbol{\mathcal{E}}$. We present
a procedure for constructing a peak function of class $\unal$. 
\medskip

The procedure now known as Bishop's one-fourth three-fourths principle \cite{bishop:mbfa59} 
(also see \cite[Theorem II/11.1]{gamelin:ua69} \ ) says that given a compact metric space
$X$, a uniform algebra $\mathcal{A}$ on $X$, and a point $x\in X$, if for each neighbourhood
$U$ of $x$, we could find a $f_U\in \mathcal{A}$ such that : {\em i)} $f_U(x)=1$, \
{\em ii)} The sup-norms $\sup_X|f_U|$ are uniformly bounded, and \ {\em iii)} 
$|f_U(y)|\leq\alpha \ \forall y\in X\setminus U$ with a uniform constant $0<\alpha<1$, then
we could construct a function $F\in \mathcal{A}$ that peaked at $x$. Our result first 
exploits the fact that the result just described can be extended beyond uniform algebras
to closed subspaces $\unal\subset_{\rm closed}\smoo$. Secondly, it develops a Bishop-type
construction in a setting where condition (ii) is replaced by a weaker analogue : 
$\sup_U|f_U|\lesssim \psi({\rm diam}(U)^{-1})$, where $\psi:\mathbb{R}_+\to\mathbb{R}_+$,
and $\psi(x)\nearrow +\infty$ sufficiently gradually as $x\to +\infty$. This sort of of
bound is motivated by applications in which the condition (ii) is difficult to verify. We
make all of this precise in the
following
\medskip

\begin{theorem}\label{T:main} Let $\Omega$ be a bounded domain in $\RN$ (or $\Cn$), and
let $x\in\OM$. Let $\unal$ be a closed subspace of $\smoo(\OM)$. Suppose there
exist constants
\begin{enumerate}
\item[{}] $0<\alpha<1$,
\item[{}] $0<s\leq 1$, $0<t<1$, and
\item[{}] $0<A<1, \ C>0$
\end{enumerate}
such that for each neighbourhood $U$ of $x$ with $\dia(U)<1$, there exists a function 
$f_U\in \unal$ with the properties
\begin{enumerate}
\item[1)] $f_U(x) \ = \ 1$;
\item[2)] $|f_U(y)| \ \leq \ \alpha \ \forall y\in \OM\setminus U$;
\item[3)] $|f_U(y)| \ \leq \ C\log^t[1/\dia(U)] \ \forall y\in U$; and
\item[4)] $\{y\in\OM \ : \ |f_U(y)|<1+\eps^s\}\supset B(x;A \ \dia(U) \ \eps), \ 0<\eps<1$.
\end{enumerate}
Then there exists a function $F\in \unal$ which peaks at $x$.
\end{theorem}
\smallskip

The expression $\dia(U)$ in the above theorem is defined as
\[
\dia (U) \ := \ \sup_{y\in U}|y-x|,
\]
where $|\centerdot |$ denotes the Euclidean norm, while the set $B(y;r)$ is defined as
$B(y;r):=\OM\cap\mathbb{B}(y;r)$, where $\mathbb{B}(y;r)$ is the open Euclidean ball 
centered at $y\in\OM$ and having radius $r$. 
\smallskip

\begin{remark} The above theorem, suitably restated, is true if $\OM$ is replaced by 
a compact metric space $X$ -- we merely state the theorem in a setting that is closest 
to its applications.
\end{remark}
\smallskip

We comment briefly on the motivation behind the precise form of condition (4) in Theorem 
\ref{T:main}. This condition is motivated by certain applications in multivariate complex 
analysis -- such as when $D$ is a pseudoconvex domain in $\CC$ having a 
sufficiently ``nice'' boundary, $x=0\in\partial{D}$, $\OM=\overline{D\cap V}$, 
$\unal=\mathcal{O}(D\cap V)\cap\smoo(\OM)$ (where $V$ is an appropriately chosen,
small $\CC$-neighbourhood of $0$), and one constructs a function of class $\unal$ that
peaks at $x=0$; refer to Fornaess \& Sibony \cite{fornaessSibony:cpfwpd89}, and
Fornaess \& McNeal \cite{fornaessMcNeal:cpfsftd94}. In both these constructions, 
the relevant $f_U$'s are scalings of a single function $f$ that satisfies (1) and (2) with 
$U=\mathbb{B}^2(0;1)\cap\overline D$ and $x=0$, and $f$ satisfies a H{\"o}lder condition 
with exponent $s=1$. The scalings $f_U$ then satisfy $|f_U(x)-f_U(y)|\lesssim \dia(U)^{-1} \ |y|$,
which {\em implies} the condition (4) with $s=1$ (but is {\em not} equivalent to (4) above).We 
ought to clarify here that the H{\"o}lder condition just stated, combined with condition (1), 
forces the relevant $f_U$'s to be uniformly bounded. Thus, in the 
constructions \cite{fornaessMcNeal:cpfsftd94} and \cite{fornaessSibony:cpfwpd89}, a 
variant of Bishop's original procedure suffices. But a different
procedure is needed in order to construct (local) peak functions on far more complicated
domains in $\Cn, \ n>2$ (such as the examples in \cite{bharali:spmdnsrt03} ) -- known as
the non-semiregular domains -- where uniform boundedness of the relevant $f_U$'s, constructed 
in analogy with \cite{fornaessSibony:cpfwpd89}, fails. However, a weaker ``H{\"o}lder-type'' 
condition, i.e. condition (4), does hold in some of these domains, whereby Theorem
\ref{T:main} can be used. Details of this last application will appear elsewhere.   
\medskip

\section{The proof of Theorem \ref{T:main}}\label{S:proofs}

To prove our theorem, we will need the following two lemmas.
\smallskip

\begin{lemma}\label{L:aux} Let $\psi : [0,\infty)\to(\eps,\infty)$ be a continuous function 
(with $\eps>0$), and assume that $\psi(x)\approx x^{1+\delta}$ for large $x$ (here, $\delta$
is a small positive number). Define
\[
g(x) \ := \ \exp\left\{-k\int_0^x\psi(s)^{-t} \ ds\right\},
\]
where $k>0$ and $t$ is such that $(1+\delta)t<1$. 
Then
\begin{enumerate}
\item[1)] There are constants $A_1,A_2>0$ such that
\begin{equation}\label{E:decay}
0 \ < \ g(x) \ \leq \ A_1\exp\left\{-\frac{kx^{1-(1+\delta)t}}{A_2(1-(1+\delta)t)}\right\}.
\end{equation}
\item[2)] $g$ satisfies the equation
\begin{equation}\label{E:inteq}
g(x) \ = \ k\int_x^{\infty}\frac{g(s)}{\psi(s)^t} \ ds.
\end{equation}
\end{enumerate}
\end{lemma}
\begin{proof}  
\medskip
Part (1) follows easily from the fact that there exist constants
$M,A_2,A_3>0$ such that 
\[
A_3^{1/t}s^{1+\delta}\ \leq \ \psi(s) \ \leq \ A_2^{1/t}s^{1+\delta} \quad 
\forall s\geq M.
\]
To prove part (2), we use the estimate \eqref{E:decay}, to see that
\[
0 < k\int_x^{\infty}\frac{g(s)}{\psi(s)^t} \ ds \ \leq \ 
\begin{cases}
{\rm const.} +
	\frac{kA_1}{A_3}\int_M^{\infty}s^{-(1+\delta)t}
	\exp\left\{-\frac{k}{A_2(1-(1-\delta)t)}s^{1-(1-\delta)t}\right\} \ ds, & \text{if 
$x<M$},\\
{} & {} \\
\frac{kA_1}{A_3}\int_x^{\infty}s^{-(1+\delta)t}
	\exp\left\{-\frac{k}{A_2(1-(1+\delta)t)}s^{1-(1-\delta)t}\right\} \ ds, 
	& \text{if $x\geq M$},
\end{cases}
\]
whence the integral on the right-hand side of \eqref{E:inteq} is a convergent integral. We make the 
following change of variable
\begin{equation}\label{E:changevar}
u(s) \ = \ k\int_0^s\psi(\tau)^{-t} \ d\tau.
\end{equation}
Since $\psi(\tau)^{-t} \ \geq \ \tau^{-(1+\delta)t}/A_2 \quad \forall \tau \geq M$, 
and $(1+\delta)t<1$, we have
\begin{equation}\label{E:atinfty}
\lim_{s\to\infty}u(s) \ = \ \infty.
\end{equation}
From \eqref{E:changevar} and \eqref{E:atinfty}, we have
\[
k\int_x^{\infty}\frac{g(s)}{\psi(s)^t} \ ds \ 	= \ \int_{u(x)}^{\infty}e^{-u} \ du \ 
						= \ e^{-u(x)} \ = \ g(x). 
\]
This completes the proof.
\end{proof}
\medskip

\begin{proof}[{\bf The proof of Theorem 1.1}] We begin by observing that we may as well
assume that $1/2<t<1$, and that $x=0$. We may henceforth assume -- by raising
the value of $C$ if necessary -- that $s-(\vaph)/C > 0$. This allows us to choose a
$D\in(0,1)$, $D$ sufficiently small, so that : 
\begin{equation}\label{E:1stShell}
\frac{\vaph}{2}D^{(\vaph )/C} \ \geq \ D^s.
\end{equation}
Recall that $A\in(0,1)$. Using this fact, we define  
\begin{align}
\log(\eps_k/D) \ &:= \ (\log{A}) \ \sum_{j=1}^{k}j^{-1+p},\quad k=1,2,\dots \notag \\
U_k \ &:= \ \begin{cases}
		B(0;D), & \text{if $k=1$},\\
		B(0;A\ \dia(U_{k-1})\ \eps_{k-1}), & \text{if $k\geq 2$},
		\end{cases} \notag
\end{align}
where we choose $p$ to satisfy $0<p<1$. Notice that $\eps_k\downarrow 0$ as $k\to\infty$.
We choose $p$ to be so close to $0$ that $(1+p)t<1$. This is always possible because
$t<1$. In fact, we shall demand that :
\begin{itemize}
\item $(1+p)t =  q$, with $q\in (0,1)$; and
\item $(1-q)=p/M$, where $M$ satisfies $M\geq C$, and whose precise value will be
stated later in this proof.
\end{itemize}
These imply that
\begin{equation}\label{E:peeque}
q \ = \ \frac{t+Mt}{1+Mt}, \qquad\quad p \ = \ \frac{M(1-t)}{1+Mt}.
\end{equation}
Observe that since we have assumed that $1/2<t<1$, the value of $p$ will indeed
be less than $1$.
\smallskip

We now construct a sequence of functions $\{f_n\}_{n\in\mathbb{N}}\subset\unal$ by 
induction as follows :
\smallskip

Define $f_1:=f_{U_1}$ and write
\[
W_1 \ = \ \{y\in \OM : |f_1(y)|\geq 1+\eps_1^s\}.
\]
Due to the condition (4) of our hypothesis (note that $\eps_1<1$), 
$W_1\subset(\OM\setminus U_2)$. If we 
define $f_2:=f_{U_2}$, it follows from the hypotheses of our theorem that :
\begin{enumerate}
\item[$(a)$] $f_2(0) \ = \ 1$;
\item[$(b)$] $\{y\in \OM : |f_2(y)|<1+\eps_2^s \} \ \supset \ U_3$;
\item[$(c)$] $|f_2(y)| \ \leq \ C\log^t[1/\dia(U_2)] \ \forall y\in U_2$; and
\item[$(d)$] $|f_2(y)| \ \leq \ \alpha \ \forall y\in (\OM\setminus U_2)$.
\end{enumerate} 
Assume that we have found $f_2,\dots,f_m\in\unal$ such that, defining
\[
W_m \ = \ \{y\in \OM : \max_{1\leq j\leq m}|f_j(y)|\geq 1+\eps_m^s\},
\]
they satisfy
\begin{enumerate}
\item[$(a)_m$] $f_j(0) \ = \ 1, \ j=2,\dots,m$;
\item[$(b)_m$] 
$\{y\in \OM : |f_j(y)|<1+\eps_j^s \} \ \supset \ U_{j+1}, \ 
j=2,\dots,m$;
\item[$(c)_m$] $|f_j(y)| \ \leq \ C\log^t[1/\dia(U_j)] \ 
\forall y\in U_j \ \text{and $j=2,\dots,m$}$; and
\item[$(d)_m$] $|f_j(y)| \ \leq \ \alpha \ \forall y\in 
(\OM\setminus U_j) \ \text{and $j=2,\dots,m$}$.
\end{enumerate} 
Notice that by $(b)_m$, $W_m\subset(\OM\setminus U_{m+1})$. If we define 
$f_{m+1}:=f_{U_{m+1}}$, then, by our hypotheses,
$\{f_2,\dots,f_m,f_{m+1}\}\subset\unal$ satisfies $(a)_{m+1}$--$(d)_{m+1}$.
\medskip

It is easy to check that for $m\geq 2$ 
\begin{align}
\log\frac{1}{\dia(U_m)} \ &= \ m\log\left(\frac{1}{D}\right) \ + \ 
		\log\left(\frac{1}{A}\right)\left\{(m-1)+
		\sum_{j=1}^{m-1}(m-j)j^{-1+p}\right\}\notag \\
&= \ m\log\left(\frac{1}{D}\right) \ + \ 
		\log\left(\frac{1}{A}\right)\left\{(m-1)+m\sum_{j=1}^{m-1}j^{-1+p}-
		\sum_{j=1}^{m-1}j^p\right\}. \label{E:radpow}
\end{align}

Notice that by estimating the sums obtained above by integrals, we have 
\begin{align}\label{E:intermest}
\frac{m^p-1}{p} \ &< \ \sum_{j=1}^{m-1}j^{-1+p} \ < \ \frac{m^p-1}{p}+(1-m^{-1+p}), \\
\frac{m^{p+1}-1}{p+1}+(1-m^p) \ &< \ \sum_{j=1}^{m-1}j^p \ < \ 
				\frac{m^{p+1}-1}{p+1},\quad m\geq 2. \notag
\end{align}
From \eqref{E:radpow} and \eqref{E:intermest}, we get
\[
\log\frac{1}{\dia(U_m)} \ \leq \ m\log\left(\frac{1}{D}\right) \ + \ 
		\log\left(\frac{1}{A}\right)\left\{2(m-1)+
		\frac{m^{1+p}-(p+1)m+p}{p(p+1)}\right\}\quad \forall m\geq 2.
\]
Therefore, there exists a constant $L>0$ that is independent of $p\in (0,1)$ such that
\begin{equation}\label{E:radpowest}
\log\frac{1}{\dia(U_m)} \ \leq \ m\log\left(\frac{1}{D}\right) \ + \ 
		\log\left(\frac{1}{A}\right)\frac{Lm^{1+p}}{p(p+1)}\quad \forall m\geq 1.
\end{equation}
\smallskip 

Define
\[
\psi(\tau) \ := \ \begin{cases}
			\tau\log(1/D) \ + \ 
			\log(1/A)\dfrac{L\tau^{1+p}}{p(p+1)}, 
				& \text{if $\tau\geq 1$}, \\
			\psi(1), & \text{if $0\leq\tau<1$}.
			\end{cases}
\]
Note that $\psi(m)\geq \log[1/\dia(U_m)] \ \forall m\in\mathbb{N}$. Finally, define
\begin{equation}\label{E:peakfunc}
F(y) \ := \ \sig^{-1}\left[\sum_{j=1}^\infty\sig_jf_j(y)\right],
\end{equation}
where
\begin{align}
g(x) \ &:= \ \exp\left\{-\frac{1-\alpha}{2M}\int_0^x\psi(s)^{-t} \ ds\right\}
\quad\text{and} \quad\sig_j \ := \frac{g(j)}{M\psi(j)^t}, \notag \\
\sig \ &:= \ \sum_{j=1}^\infty\sig_j. \notag
\end{align}
It is easy to check that the last series above is rapidly convergent. To see this, we
apply Lemma \ref{L:aux} to $\psi$. By the manner in which $p$ and $q$ are defined, we have the
estimate
\[
0 \ < \ \sig_j \ \leq \ 
\frac{A_1}{M\psi(j)^t}\exp\left\{-\frac{1-\alpha}{2MA_2(1-q)} \ j^{1-q}\right\},
\]
where $A_1$ and $A_2$ are the constants given by Lemma \ref{L:aux}. Thus, by item (3) of
our hypothesis, 
\[
0 \ < \ \sigma_j\sup_{\OM}|f_j| \ \leq \ 
	\frac{CA_1}{M} \ \exp\left\{-\frac{1-\alpha}{2MA_2(1-q)} \ j^{1-q}\right\}.
\]
Since the right-hand side of the above estimate constitutes a summable series, as $j$
varies over $\mathbb{N}$, we conclude that the right-hand side of \eqref{E:peakfunc}
converges uniformly on $\OM$. Therefore $F\in\unal$.
\medskip
     
We claim that $F$ defined by \eqref{E:peakfunc} peaks at $x$. Before proving this 
assertion, we choose an appropriately large value for $M$, which links $p$ and $q$ via
the relation
\[
(1-q)=p/M.
\]
We choose $M$ to be so large -- i.e. $q<1$ to be so close to $1$ -- that :
\begin{align}
M &\geq C, \quad and \notag \\
\exp\left\{-\frac{1+Mt}{1-t} \ \right.&
		\left.\frac{\vaph}{2M}[ \ (m+1)^{p/M}-2^{p/M}]\right\} \label{E:eps} \\ 
	&\geq \ \exp\left\{-\log\frac{1}{A} \ \frac{s(1+Mt)}{M(1-t)}
			[ \ (m-1)^p-1 \ ]\right\}\quad \forall m\geq 3. \notag
\end{align}
With this choice of $p$ and $q$, we can make the following claims :
\smallskip

\noindent{{\bf Claim 1.} $(C\psi(m)^t-1)\sig_m < (\vaph)\sum_{j\geq m+1}\sig_j/2$ for
each $m\in\mathbb{N}$.}
\medskip

\noindent{{\bf Claim 2.} $(\vaph)\sum_{j\geq m+1}\sig_j/2 > \eps_{m-1}^s\sum_{j=1}^{m-1}\sig_j$
for each $m\geq 2$.}
\medskip

We defer the proofs of these claims to the end of this section. Assuming that these claims are
true, we can show that $F$ peaks at $0$. We first consider $y\neq 0$ and 
$y\notin\cup_{j\in\mathbb{N}}W_j$. Then, $|f_j(y)|\leq 1 \ \forall j\in\mathbb{N}$.
However, $y\in(\OM\setminus U_{j_0})$ for some $j_0\in\mathbb{N}$,
whence, by condition $(d)_m$ : $|f_{j_0}(y)|\leq\alpha<1$. Consequently, $|F(y)|<1$. This
leaves us with the case $y\in\cup_{j\in\mathbb{N}}$ to analyze. Since $\{W_j\}_{j\in\mathbb{N}}$
is a strictly increasing sequence of
closed sets, either $y\in W_1$ or there exists $m\in\mathbb{N}$ such that
$y\in W_m$ but $y\notin W_j \ \forall j\leq m-1$. In the former case, it is clear by
construction that $|F(y)|\leq\alpha<1$. The latter case results in the following estimates
\begin{align}
|f_j(y)| \ &< \ 1+\eps_{m-1}^s, \quad 1\leq j<m, \notag \\
|f_m(y)| \ &\leq \ C\log^t[1/\dia(U_m)] \ = \ C\psi(m)^t, \notag \\
|f_j(y)| \ &\leq \ \alpha \quad \forall j\geq m+1. \notag
\end{align}
Then
\begin{align}
|F(y)| \ 
&\leq \ \sig^{-1}\left\{(1+\eps_{m-1}^s)\sum_{j=1}^{m-1}\sig_j+C\psi(m)^t\sig_m
		+\alpha\sum_{j=m+1}^\infty\sig_j\right\}	&& {} \notag \\
&< \ \sig^{-1}\left\{\left(\sum_{j=1}^{m-1}\sig_j+
	\frac{1-\alpha}{2}\sum_{j=m+1}^\infty\sig_j\right)\right.	&& {} \notag \\
&\qquad\qquad\quad\left.+\left(\sig_m+\frac{1-\alpha}{2}\sum_{j=m+1}^\infty\sig_j\right)
		+\alpha\sum_{j=m+1}^\infty\sig_j\right\} \ = \ 1.	
	&&\text{(from Claims 1 \& and 2)} \notag 
\end{align}
Thus, $|F(y)|<1 \ \forall y\neq 0$, whence $F$ peaks at $0$.
\medskip

To complete our proof, we first present
\smallskip

\noindent{{\bf The proof of Claim 1 :} We compute :
\begin{align}
(C\psi(m)^t-1)\sig_m \ &= \ \frac{C}{M} \ g(m)-\frac{g(m)}{M\psi(m)^t} 	&& {} \notag \\
&< \ \frac{C(1-\alpha)}{2M}\int_m^{\infty}\frac{g(s)}{M\psi(s)^t} \ ds
	- \frac{1-\alpha}{2}\sig_m 	&&\text{(applying Lemma \ref{L:aux}-2)} \notag \\
&< \ \frac{1-\alpha}{2}\sum_{j=m}^{\infty}\frac{g(j)}{M\psi(j)^t}
	- \frac{1-\alpha}{2}\sig_m	&&\text{(estimating from above by a series)} \notag \\
&= \ \frac{1-\alpha}{2}\sum_{j=m+1}^{\infty}\sig_j.	&& {} \notag
\end{align}
This proves Claim 1.}
\medskip

And finally, we present
\smallskip

\noindent{{\bf The proof of Claim 2 :} Notice that when $m\geq 2$,
\begin{align}
\frac{\frac{1-\alpha}{2}\sum_{j=m+1}^\infty\sig_j}{\sum_{j=1}^{m-1}\sig_j} \ 
&> \ \frac{\frac{1-\alpha}{2}\int_{m+1}^\infty g(s)/M\psi(j)^t \ ds}
	{\int_0^{m-1} g(s)/M\psi(j)^t \ ds}	&&\text{(estimating sums by integrals)} 
\notag \\
&> \ \frac{1-\alpha}{2} \ \frac{g(m+1)}{g(0)-g(m-1)}	&&\text{(applying Lemma \ref{L:aux}-2)} \notag \\
&> \ \frac{1-\alpha}{2} \ \frac{g(m+1)}{g(0)} 	&& {} \notag \\
&= \ \frac{1-\alpha}{2} \ \exp\left\{-\frac{1-\alpha}{2M}\int_0^{m+1}\psi(s)^{-t} \ ds\right\}
	&& {} \notag \\
&\geq \frac{1-\alpha}{2} \ D^{(\vaph)/M}
		\exp\left\{-\frac{1-\alpha}{2M}\int_2^{m+1}s^{-t(1+p)} \ ds\right\}
	&&\text{(since $s^{1+p}\leq \psi(s) \ \forall s\geq 2$)} \notag \\
&\geq \ \frac{1-\alpha}{2} \ D^{(\vaph)/C} && {} \label{E:claim2} \\
&\qquad\qquad\times\exp\left\{-\frac{1-\alpha}{2M} \ \frac{(m+1)^{1-q}-2^{1-q}}{1-q}\right\}
	&&\text{(since $C\leq M$ and $0<D<1$)} \notag
\end{align}
At this stage, we use the condition \eqref{E:eps} (recall that $m\geq 2$) to get
\begin{align}\label{E:last}
\frac{1-\alpha}{2} & D^{(\vaph)/C}
\exp\left\{-\frac{1-\alpha}{2M} \ \frac{(m+1)^{1-q}-2^{1-q}}{1-q}\right\} && {} \\
&= \ \frac{1-\alpha}{2} \ D^{(\vaph)/C} && {} \notag \\
& \ \qquad\times
\exp\left\{-\frac{1+Mt}{1-t} \ \frac{\vaph}{2M} \ [ \ (m+1)^{p/M}-2^{p/M}]\right\} 
&& {} \notag \\
&\geq \ D^s\exp\left\{-\log\frac{1}{A} \ \frac{s}{p} \ 
			[ \ (m-1)^p-1 \ ]\right\} 
			&&\text{(using the fact $(1+Mt)/M(1-t)=1/p$)} \notag \\
&\geq \ D^sA^{s\sum_{j\leq(m-1)}j^{-1+p}} 
		&&\text{(using \eqref{E:intermest} on the exponents)} \notag \\
&= \ \eps_{m-1}^s. && {} \notag
\end{align}
Note that the first inequality makes use of the condition \eqref{E:1stShell}.
Comparing \eqref{E:last} with \eqref{E:claim2} we conclude that
\[
\frac{1-\alpha}{2}\sum_{j=m+1}^\infty\sig_j \ > \ \eps_{m-1}^s\sum_{j=1}^{m-1}\sig_j\quad
\forall m\geq 2,
\]
which is precisely Claim 2.}
\medskip

This concludes our proof.
\medskip
\end{proof}

\end{document}